\documentclass[12pt]{article}
\usepackage{amsmath}
\usepackage{amsfonts}
\usepackage{mathrsfs}
\usepackage{amssymb}
\usepackage{indentfirst}
\usepackage[usenames]{color}
\usepackage[numbers,sort&compress]{natbib}
\usepackage{subfigure}
\usepackage{bm}
\usepackage{multirow}
\usepackage{setspace}
\usepackage[all]{xy}

\def\cl{\centerline}
\def\ni{\noindent}

\def\BB{{\cal B}(q)}
\def\sp{{\rm{span}}}
\def\deg{{\rm deg}}

\def\Z{\mathbb{Z}}

\def\C{\mathbb{C}}

\def\QED{\hfill$\Box$}

\def\Uh0{\mathcal{U}(\mathfrak h_0)}
\def\Bp0{\mathcal{B}(p,0)}
\def\bp1{\mathcal{B}(p,1)}

\numberwithin{equation}{section}
\newtheorem{theo}{Theorem}[section]

\newtheorem{lemm}[theo]{Lemma}
\newtheorem{rema}[theo]{Remark}
\newtheorem{prop}[theo]{Proposition}

\newtheorem{conv}[theo]{Convention}

\begin{document}
\begin{center}
{\bf\large Representations of non-finitely graded Heisenberg-Virasoro type Lie algebras}
\footnote{
$^{*}$Corresponding author: chgxia@cumt.edu.cn (C.~Xia).}
\end{center}
\vspace{8pt}

\cl{\footnotesize Chunguang Xia$^{\,*,\,\dag,\,\ddag}$, Tianyu Ma$^{\,\dag}$, Wei Wang$^{\,\S}$, Mingjing Zhang$^{\,\dag,\,\ddag}$}
\vspace{8pt}

\cl{\footnotesize $^{\,\dag}$School of Mathematics, China University of Mining and Technology, Xuzhou 221116, China}
\vspace{8pt}
\cl{\footnotesize $^{\,\ddag}$Jiangsu Center for Applied Mathematics (CUMT), Xuzhou 221116, China}
\vspace{8pt}
\cl{\footnotesize $^{\,\S}$School of Mathematics and Statistics, Ningxia University, Yinchuan 750021, China}
\vspace{8pt}
\cl{\footnotesize E-mails: chgxia@cumt.edu.cn, tyma@cumt.edu.cn,   wwll@mail.ustc.edu.cn, zhangmj@nju.edu.cn}
\vspace{10pt}

\footnotesize
\noindent{{\bf Abstract.}
We construct and study non-finitely graded Lie algebras $\mathcal{HV}(a,b;\epsilon)$ related to Heisenberg-Virasoro type Lie algebras,
where $a,b$ are complex numbers, and $\epsilon = \pm 1$.
Using combinatorial techniques, we completely classify the free $\mathcal{U}(\mathfrak h)$-modules of rank one over $\mathcal{HV}(a,b;\epsilon)$.
It turns out that these modules are more varied and complex than those over non-finitely graded Virasoro algebras, and in particular
admit infinitely many free parameters if $b=1$ and $\epsilon=-1$.
Meanwhile, we also determine the simplicity and isomorphism classes of these modules.

\ni{\bf Key words:} Heisenberg-Virasoro algebra; non-finitely graded Lie algebras; free $\mathcal{U}(\mathfrak h)$-module; simplicity

\ni{\it Mathematics Subject Classification (2020):} 17B10; 17B65; 17B68}

%\small
%\begin{spacing}{1}
%\renewcommand{\baselinestretch}{1}
%\tableofcontents
%\end{spacing}

\small
\section{\large{Introduction}}

The Virasoro algebra is a $\Z$-graded Lie algebra, which plays an important role in many areas of mathematics and physics.
%which has historically appeared in both mathematical and physical context, and plays a central role in conformal field theory, especially during 1980s.
Its representation theory has been extensively and deeply studied, especially Harish-Chandra modules \cite{M1992} and Verma modules \cite{FF1990}.
Following the spirit of \cite{FF1990,M1992}, one can logically study the representations of finitely graded Lie algebras with Virasoro subalgebra,
such as those of (twisted, generalized) Heisenberg-Virasoro Lie algebras (e.g.,~\cite{ADKP1998,B2003,SJS2008,LG2014,AR2015}).
%and (generalized) Schr\"{o}dinger-Virasoro Lie algebras \cite{TZ2009,L2016}.

However, the representation theory of non-finitely graded Lie algebras with Virasoro subalgebra is rather difficult in general.
The first important breakthrough in this direction is to study the quasifinite modules, 
which was initiated by Kac and Radul \cite{KR1993} for the $W$-infinity algebra $\mathcal{W}_{1+\infty}$,
see also subsequent work for Lie algebras of Weyl type \cite{S2003} and Lie algebras of Block type \cite{WT2007,SXX2012,SXX2013}.
A common feature of these Lie algebras is that their zero-graded parts with respect to principal gradations are commutative subalgebras.

For non-finitely graded Lie algebras with no the above feature, the new difficulty is that one cannot cope with quasifinite modules as usual.
The prototypical example is the so-called non-finitely graded Virasoro algebra $\mathcal{W}(\epsilon)$, $\epsilon=\pm1$ \cite{SXZ2000,GJ2010,XMDZ2023}.
To overcome this difficulty, one can turn to consider the modules satisfying certain ``generalized quasifinite conditions'' \cite{GJ2010}.
Unfortunately, beside \cite{GJ2010}, there is no further research on the representations on this kind of Lie algebras during the last fifteen years.
A new way we recently found in \cite{XMDZ2023} is to study the free $\mathcal{U}(\mathfrak h)$-modules.
This kind of modules \cite{N2015,N2016} was originally defined by an ``opposite condition'' relative to weight modules,
and constructed first for $\mathfrak{sl}_n$ by Nilsson \cite{N2015} and independently by Tan and Zhao \cite{TZ2018}.
Usually, this kind of modules admits finitely many free parameters, see, e.g., \cite{CC2015,N2016,TZ2015,XMDZ2023} and the references therein.

Let $a,b\in\C$ and $\epsilon=\pm 1$. In this paper,
we construct and study a class of non-finitely graded Lie algebras $\mathcal{HV}(a,b;\epsilon)$,
which are referred to as {\it non-finitely graded Heisenberg-Virasoro type Lie algebras},
considering that each $\mathcal{HV}(a,b;\epsilon)$ contains two important subalgebras. One is non-finitely graded Virasoro subalgebra
$\mathcal{W}(\epsilon)$ mentioned above (the prototypical example of non-finitely graded Lie algebras which can not be studied through quasifinite modules), 
and the other is Heisenberg-Virasoro type subalgebra $\mathfrak{hv}(a,b)$.
Also, $\mathcal{W}(\epsilon)$ and $\mathfrak{hv}(a,b)$ share a centerless Virasoro subalgebra $\mathfrak{vir}$.
Namely, we have the following inclusion relations:
$$
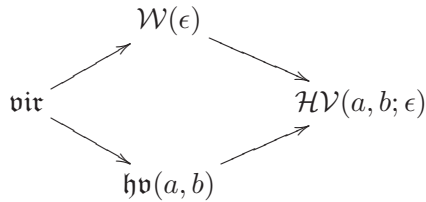
\xymatrix@R=0.5cm{
                &         \mathcal{W}(\epsilon)\ar[dr]   &   \\
\mathfrak{vir} \ar[ur] \ar[dr]           & &   \mathcal{HV}(a,b;\epsilon) \\
                &        \mathfrak{hv}(a,b)\ar[ur]  &              }
$$
See Subsection~2.1 for more details.
Due to the importance of the above mentioned subalgebras, the present study was conducted to understand the representations of $\mathcal{HV}(a,b;\epsilon)$.

Take $\mathfrak{h}=\C L_{0,0}\in\mathcal{HV}(a,b;\epsilon)$. Using combinatorial techniques,
we completely classify the free $\mathcal{U}(\mathfrak h)$-modules of rank one over $\mathcal{HV}(a,b;\epsilon)$.
In \cite{XMDZ2023}, we found that the module structures of $\mathcal{W}(\epsilon)$ are quite different from those of known non-finitely graded Lie algebras.
In this paper, we show that the module structures of $\mathcal{HV}(a,b;\epsilon)$ are more varied and complex (Theorems~\ref{thm-hvab1} and \ref{thm-hvab-1}).
In particular, there will appear modules admitting infinitely many free parameters if $b=1$ and $\epsilon=-1$.
 This is a {\bf rather unusual phenomenon} in contrast to previous results.
Denote by $\Omega_\mathfrak{g}$ the free $\mathcal{U}(\mathfrak h)$-modules of rank one over
the Lie algebra $\mathfrak{g}=\mathfrak{vir},\mathfrak{hv}(a,b),\mathcal{W}(\epsilon),\mathcal{HV}(a,b;\epsilon)$,
and by $\sharp(\Omega_\mathfrak{g})$ the number of free parameters of $\Omega_\mathfrak{g}$.
Comparing our results (and recalling Convention~\ref{conv-iso}) with those in \cite{TZ2015,CC2015,XMDZ2023}, we have Tables 1 and 2.

\begin{table}[h]\label{tab-1}
\centering\small
\subtable[$\sharp(\Omega_\mathfrak{g})$ for $\mathfrak{g}=\mathfrak{vir}$ and $\mathcal{W}(\epsilon)$]{
\begin{tabular}{c|c|c}
\hline{\centering} \rule{0pt}{11pt}
Lie algebra $\mathfrak{g}$  & $\sharp(\Omega_\mathfrak{g})$ &  Reference \\
\hline \rule{0pt}{11pt}
$\mathfrak{vir}$ &  2 &  \cite[Theorem~3]{TZ2015} \\
\hline \rule{0pt}{11pt}
$\mathcal{W}(\epsilon)$ & 3 &  \cite[Theorems~3.2 and 4.2]{XMDZ2023} \\
\hline
\end{tabular}}
\end{table}

\begin{table}[h]\label{tab-2}
\centering\small
\subtable[$\sharp(\Omega_\mathfrak{g})$ for $\mathfrak{g}=\mathfrak{hv}(a,b)$ and $\mathcal{HV}(a,b;\epsilon)$]{
\begin{tabular}{c|c|c|c|c|c}
\hline{\centering} \rule{0pt}{11pt}
Lie algebra $\mathfrak{g}$  & $a$ &  $b$ & $\epsilon$ & $\sharp(\Omega_\mathfrak{g})$ & Reference \\
\hline \rule{0pt}{11pt}
\multirow{4}*{$\mathfrak{hv}(a,b)$ } & $a$ &  $b=0$ & \multirow{4}*{--} &  3 & \multirow{4}*{\cite[Theorem~2.5]{CC2015}}\\
\cline{2-3} \rule{0pt}{11pt}
~& $a\ne0$ &   $b=1$ &  & 3 & ~\\
\cline{2-3} \rule{0pt}{11pt}
~& $a=0$ &   $b=1$ &  & 3 & ~\\
\cline{2-3} \rule{0pt}{11pt}
~& $a$ &   $b\ne 0,1$ &  & 2 & ~\\
\hline \rule{0pt}{11pt}
\multirow{9}*{$\mathcal{HV}(a,b;\epsilon)$} & $a$ &  $b=0$ & \multirow{4}*{1} &  4 & \multirow{4}*{Theorem~\ref{thm-hvab1}}\\
\cline{2-3} \rule{0pt}{11pt}
~& $a\ne0$ &   $b=1$ &  & 4 & ~\\
\cline{2-3} \rule{0pt}{11pt}
~& $a=0$ &   $b=1$ &  & 3 & ~\\
\cline{2-3} \rule{0pt}{11pt}
~& $a$ &   $b\ne 0,1$ &  & 3 & ~\\
\cline{2-6} \rule{0pt}{11pt}
~ & $a$ &  $b=0$ & \multirow{4}*{$-1$} &  4 & \multirow{4}*{Theorem~\ref{thm-hvab-1}}\\
\cline{2-3} \rule{0pt}{11pt}
~& $a\ne0$ &   $b=1$ &  & $\infty$ & ~\\
\cline{2-3} \rule{0pt}{11pt}
~& $a=0$ &   $b=1$ &  & $\infty$ & ~\\
\cline{2-3} \rule{0pt}{11pt}
~& $a$ &   $b\ne 0,1$ &  & 3 & ~\\
\hline
\end{tabular}}
\end{table}

Our proofs in this paper and \cite{XMDZ2023} indicate that one need to employ various combinatorial techniques
when dealing with representations of Lie algebras related to $\mathcal{W}(\epsilon)$, such as Pascal's triangle.
Moreover, in order to present the classification results in this paper in a more compact form, we utilize generalized binomial coefficients
(together with higher-order derivatives). Thus one may expect more combinatorial ingredients to appear in the future study related to this kind of Lie algebras. This is also our motivation for presenting this work.

This paper is organized as follows. In Section 2, we construct the Lie algebra $\mathcal{HV}(a,b;\epsilon)$,
and recall the results of free $\mathcal{U}(\mathfrak h)$-modules over Lie algebras $\mathcal{W}(\epsilon)$ and $\mathfrak{hv}(a,b)$.
In Sections 3 and 4, by using combinatorial techniques, we completely classify the free $\mathcal{U}(\mathfrak h)$-modules of rank one over $\mathcal{HV}(a,b;\epsilon)$
(Theorems~\ref{thm-hvab1} and \ref{thm-hvab-1}).
Meanwhile, we also determine the simplicity and isomorphism classes of these modules (Theorems~\ref{hvab1-sim&iso} and \ref{hvab-1-sim&iso}).

\section{\large{Preliminaries}}

Throughout this paper, we use $\Z$, $\Z_+$, $\C$ and $\C^*$ to
denote the sets of integers, nonnegative integers, complex numbers and nonzero complex numbers, respectively.
In this section, we first construct the Lie algebra $\mathcal{HV}(a,b;\epsilon)$,
and then recall the free $\mathcal{U}(\mathfrak h)$-modules results of $\mathcal{W}(\epsilon)$ \cite{XMDZ2023} in a unified and compact form.
We also recall the free $\mathcal{U}(\mathfrak h)$-modules results of $\mathfrak{hv}(a,b)$ \cite{CC2015} for later use.

\subsection{\normalsize{Construction of the Lie algebra $\mathcal{HV}(a,b;\epsilon)$}}

Let $\epsilon=\pm 1$. We first recall two kinds of realizations of the non-finitely graded Virasoro algebra $\mathcal{W}(\epsilon)$ \cite{CHS2014,XMDZ2023}.
One is an algebraic realization and the other one is a geometric realization.

(1) {\bf Algebraic realization.} Consider the semigroup algebra $\C[x^{\pm1},t]$. Let $\partial_x$ and $\partial_t$ be the usual partial derivatives of $\C[x^{\pm1},t]$.
Denote by $\partial_\epsilon = x^2\partial_x + t^{1-\epsilon}\partial_t$. Then, one can check that $\mathcal{W}(\epsilon)$ is the following Lie algebra
$$
\mathcal{W}(\epsilon)=\C[x^{\pm1},t]\partial_\epsilon=\sp\{L_{i,m} := -\epsilon x^{-\epsilon i}t^m\partial_\epsilon |\, i\in\Z,m\in\Z_+\}.
$$

(2) {\bf Geometric realization.} We know that the space $C^{\infty}_{[0,+\infty)}$ of smooth function in the  interval $[0,+\infty)$
becomes a Lie algebra under bracket $[f,g]=fg'-f'g$. One can check that $\mathcal{W}(\epsilon)$ is exactly the following subalgebra of $C^{\infty}_{[0,+\infty)}$:
$$
\mathcal{W}(\epsilon)=\sp\{L_{i,m}:=-(1+t)^{\epsilon m} e^{-it}\in C^{\infty}_{[0,+\infty)}\, |\, i\in\Z,m\in\Z_+\}.
$$

Let $a,b\in\C$. It is straightforward to check that the space $V_{a,b}=\sp\{v_{j,n}\mid j\in\Z,n\in\Z_+\}$
becomes a $\mathcal{W}(\epsilon)$-module under the action (see \cite{SZ2002,GJ2010,WY2021} for some special cases)
\begin{equation*}
L_{i,m}v_{j,n} = (a + j + bi)v_{i+j,m+n} - \epsilon (n+bm)v_{i+j,m+n-\epsilon}.
\end{equation*}
Motivated by this module structure, it is reasonable to consider the space
$$
\mathcal{HV}(a,b;\epsilon)=\sp\{L_{i,m},H_{j,n}\, |\,i,j\in\Z, m,n\in\Z_+\},
$$
satisfying relations
\begin{eqnarray}
\label{hv-def-1}  [L_{i,m}, L_{j,n}] &\!\!\!=\!\!\!& (j-i)L_{i+j,m+n} - \epsilon(n-m)L_{i+j,m+n-\epsilon}, \\
\label{hv-def-2}  [L_{i,m}, H_{j,n}] &\!\!\!=\!\!\!& (a+j+bi)H_{i+j,m+n} - \epsilon(n+bm)H_{i+j,m+n-\epsilon}, \\
\label{hv-def-3}  [H_{i,m}, H_{j,n}] &\!\!\!=\!\!\!& 0.
\end{eqnarray}
Obviously, it is a Lie algebra isomorphic to the semidirect product $\mathcal{W}(\epsilon)\ltimes V_{a,b}$.
Note that it contains a Heisenberg-Virasoro type subalgebra $\mathfrak{hv}(a,b)$ spanned by $\{L_{i,0},H_{j,0}\, |\,i,j\in\Z\}$
(the special case $\mathfrak{hv}(0,-1)$ is exactly the $W$-algebra $W(2,2)$ \cite{ZD2009},
and the special case $\mathfrak{hv}(0,0)$ is the centerless twisted Heisenberg-Virasoro Lie algebra \cite{ADKP1998}).
Hence, we refer to $\mathcal{HV}(a,b;\epsilon)$ as {\it non-finitely graded Heisenberg-Virasoro type Lie algebras}.
Note also that $\mathcal{W}(\epsilon)$ and $\mathfrak{hv}(a,b)$ share a centerless Virasoro subalgebra $\mathfrak{vir}$ spanned by $\{L_{i,0}\, |\,i\in\Z\}$.
The structures, including derivations, automorphisms, and the second cohomology group, of $\mathcal{HV}(0,0;1)$ have been studied in \cite{FZY2015}.

For simplicity, we shall adopt the following convention:

\begin{conv}\label{conv-iso}
Let $k\in\Z$.  One easily sees that the map $\Phi$ defined by
$$\Phi: \mathcal{HV}(a,b;\epsilon)\rightarrow \mathcal{HV}(a+k,b;\epsilon), \quad L_{i,m} \mapsto L'_{i,m},\quad H_{j,n} \mapsto H'_{j-k,n}$$
is a Lie algebra isomorphism.
In view of this, in what follows, we always assume that $0 \leq Re(a)< 1$, where $Re(a)$ denotes the real part of $a\in\C$.
\end{conv}

By \eqref{hv-def-1}--\eqref{hv-def-3}, it is not difficult to give a generating set of $\mathcal{HV}(a,b;\epsilon)$:
\begin{equation}\label{gen-set}
\mathcal{HV}(a,b;\epsilon)=
\left\{
\begin{aligned}
&\langle L_{i,0},\,L_{i,1},\,H_{i,0}\mid i\in\Z\rangle, && \mbox{if\ } \epsilon=1 \mbox{\ or\ } \epsilon=-1, b\ne 1,\\
&\langle L_{i,0},\,L_{i,1},\,H_{i,0},\,H_{i,1}\mid i\in\Z\rangle, && \mbox{if\ } \epsilon=-1, b=1.
\end{aligned}\right.
\end{equation}
The above difference on generating sets inspires us to ask: {\bf Whether} there are distinct module structures over $\mathcal{HV}(a,1;-1)$?
We shall give an affirmative answer in the final section of this paper.

\subsection{\normalsize{Free $\mathcal{U}(\mathfrak h)$-modules of rank one over $\mathcal{W}(\epsilon)$}}

Let $m\in\Z$ and $n\in\Z_+$. The {\it generalized binomial coefficient} $\binom{m}{n}$ is defined by
$$
\binom{m}{n}=
\left\{
\begin{aligned}
&1, && \mbox{if\ } n=0,\\
&\frac{1}{n!}\prod_{i=m-n+1}^{m}i, && \mbox{if\ } n>0.
\end{aligned}\right.
$$
Note that, if $m>0$, then $\binom{-m}{n}=(-1)^n\binom{m+n-1}{n}$.
For any polynomial $f(t)\in\C[t]$, we shall use $f^{(n)}(t)$ to denote the $n$-order derivative of $f(t)$.
Recall that the Lie algebra $\mathcal{W}(\epsilon)$ has brackets \eqref{hv-def-1}, and $\mathfrak h=\C L_{0,0}\in\mathcal{W}(\epsilon)$.
Using combinatorial techniques, we have shown in \cite{XMDZ2023} that

\begin{lemm}\label{We-module}
Any free $\mathcal{U}(\mathfrak h)$-module of rank one over $\mathcal{W}(\epsilon)$ is isomorphic to
$
\Omega_{\mathcal{W}(\epsilon)}(\lambda, \alpha, \beta):=\C[t]
$
with actions
\begin{equation*}\label{action-of-We}
L_{i,m}\cdot f(t)=
\sum_{s=0}^{+\infty}\binom{\epsilon m}{s}\lambda^i\beta^{m-\epsilon s}(t-i\alpha+(\epsilon m-s)\alpha\beta^{-\epsilon} )f^{(s)}(t-i),
\end{equation*}
for some $\lambda \in\C^*$ and $\alpha, \beta\in\C$.
\end{lemm}

%\sum_{s=0}^{{\rm min}\{m,\,k\}}s!\binom{m}{s}\binom{k}{s}\lambda^i\beta^{m-s-1}((m-s)\alpha - i\alpha\beta + \beta t)(t-i)^{k-s}, && \mbox{if\ } \epsilon=1,\\
%\sum_{s=0}^{k}(-1)^ss!\binom{m+s-1}{s}\binom{k}{s}\lambda^i\beta^{m+s}(t-i\alpha-(m+s)\alpha\beta)(t-i)^{k-s}, && \mbox{if\ } \epsilon=-1

Here, we have written the $\mathcal{W}(\epsilon)$-module structure in a unified and compact form.
We refer the reader to \cite[Theorems~3.2 and 4.2]{XMDZ2023} for the original forms (see also \eqref{action-hvab1-L} for $\epsilon=1$ and \eqref{action-hvab-1-L} for $\epsilon=-1$).
Clearly, $\Omega_{\mathcal{W}(\epsilon)}(\lambda, \alpha, \beta)$ has three free parameters (as listed in Table~\ref{tab-1}).
We have also given the simplicity and isomorphism classification of $\Omega_{\mathcal{W}(\epsilon)}(\lambda, \alpha, \beta)$ in \cite[Theorems~3.5, 3.6, 4.6 and 4.7]{XMDZ2023}.

\begin{lemm}\label{we-sim&iso}The following statements hold:
\begin{itemize}\parskip-3pt
\item[{\rm(1)}] $\Omega_{\mathcal{W}(\epsilon)}(\lambda,\alpha,\beta)$ is simple if and only if $\alpha\neq0$;
\item[{\rm(2)}] $\Omega_{\mathcal{W}(\epsilon)}(\lambda,0,\beta)$ has a unique proper submodule $t\Omega_{\mathcal{W}(\epsilon)}(\lambda,0,\beta)\cong\Omega_{\mathcal{W}(\epsilon)}(\lambda,1,\beta)$, and the quotient $\Omega_{\mathcal{W}(\epsilon)}(\lambda,0,\beta)/t\Omega_{\mathcal{W}(\epsilon)}(\lambda,0,\beta)$ is a one-dimensional trivial $\mathcal{W}(\epsilon)$-module;
\item[{\rm(3)}] $\Omega_{\mathcal{W}(\epsilon)}(\lambda_1,\alpha_1,\beta_1)\cong\Omega_{\mathcal{W}(\epsilon)}(\lambda_2,\alpha_2,\beta_2)$
if and only if they have the same parameters.
\end{itemize}
\end{lemm}

\subsection{\normalsize{Free $\mathcal{U}(\mathfrak h)$-modules of rank one over $\mathfrak{hv}(a,b)$}}

Recall that the Lie algebra $\mathfrak{hv}(a,b)=\sp\{L_i:=L_{i,0},\,H_j:=H_{j,0}\mid i,j\in\Z\}$ has brackets
$$
[L_i, L_j]=(j-i)L_{i+j},\quad [L_i, H_j]=(a+j+bi)H_{i+j},\quad [H_i, H_j]=0,
$$
and $\mathfrak h=\C L_{0}\in\mathfrak{hv}(a,b)$.
Let $k\in\Z$ and recall Convention~\ref{conv-iso}. The restriction of $\Phi$ on $\mathfrak{hv}(a,b)\subset\mathcal{HV}(a,b;\epsilon)$ is also a Lie algebra isomorphism
from $\mathfrak{hv}(a,b)$ to $\mathfrak{hv}(a+k,b)$. Thus it is also reasonable to assume here that $0 \leq Re(a)< 1$.
It was shown in \cite[Theorem~2.5]{CC2015} that

\begin{lemm}\label{hvab-module}
Any free $\mathcal{U}(\mathfrak h)$-module of rank one over $\mathfrak{hv}(a,b)$ is isomorphic to
$
\Omega_{\mathfrak{hv}(a,b)}(\lambda,\alpha,\gamma):=\C[t]
$
with actions
\begin{equation*}
\begin{aligned}\label{action-of-hvab}
L_i\cdot f(t) &= \lambda^i(t-i\alpha)f(t-i),\\
H_i\cdot f(t) &= \phi_i(a,b)\lambda^i\gamma f(t-a-i),
\end{aligned}
\end{equation*}
for some $\lambda\in\C^*$ and $\alpha, \gamma\in\C$, where
\begin{equation}\label{func-phiiab}
\phi_i(a,b)=
\left\{\begin{array}{ll}
1, & b=0,\\[2pt]
\frac{a}{a+i}, & b=1,\, a \neq 0,\\[4pt]
\delta_{i,0}, & b=1,\, a = 0,\\[4pt]
0, & b\neq 0,1.\\[4pt]
\end{array}\right.
\end{equation}
\end{lemm}

Note that $\Omega_{\mathfrak{hv}(a,b)}(\lambda,\alpha,\gamma)$ has three free parameters if $b=0,1$, and two free parameters if $b\ne0,1$ (as listed in Table~2).
The simplicity and isomorphism classification of $\Omega_{\mathfrak{hv}(a,b)}(\lambda,\alpha,\gamma)$ are as follows \cite[Propositions~2.3 and 2.4]{CC2015}.

\begin{lemm}\label{hvab-sim&iso}The following statements hold:
\begin{itemize}\parskip-3pt
\item[{\rm(1)}] $\Omega_{\mathfrak{hv}(a,b)}(\lambda,\alpha,\gamma)$ is simple if and only if $\alpha \neq 0\ {\rm or}\ \gamma\neq 0$ (if $\gamma$ exists);
\item[{\rm(2)}] $\Omega_{\mathfrak{hv}(a,b)}(\lambda,0,0)$ has a unique proper submodule $t\Omega_{\mathfrak{hv}(a,b)}(\lambda,0,0)\cong\Omega_{\mathfrak{hv}(a,b)}(\lambda,1,0)$, and the quotient $\Omega_{\mathfrak{hv}(a,b)}(\lambda,0,0)/t\Omega_{\mathfrak{hv}(a,b)}(\lambda,0,0)$ is a one-dimensional trivial $\mathfrak{hv}(a,b)$-module;
\item[{\rm(3)}] $\Omega_{\mathfrak{hv}(a,b)}(\lambda_1,\alpha_1,\gamma_1)\cong\Omega_{\mathfrak{hv}(a,b)}(\lambda_2,\alpha_2,\gamma_2)$ if and only if they have the same parameters.
\end{itemize}
\end{lemm}

\section{\large{Free $\mathcal{U}(\mathfrak h)$-modules of rank one over $\mathcal{HV}(a,b;1)$}}

Recall that $\mathfrak h=\C L_{0,0}\in \mathcal{HV}(a,b;1)$.
In this section, we completely classify the free $\mathcal{U}(\mathfrak h)$-modules of rank one over $\mathcal{HV}(a,b;1)$.
One will see that the module structures of $\mathcal{HV}(a,b;1)$ are more varied and complex than those of $\mathcal{W}(1)$.

\subsection{\normalsize{Construction of free $\mathcal{U}(\mathfrak h)$-modules of rank one over $\mathcal{HV}(a,b;1)$}}

Let $\lambda\in\C^*$ and $\alpha,\beta,\gamma\in\C$.
Define the action of $\mathcal{HV}(a,b;1)$ on the vector space of polynomials in one variable
$
\Omega_{\mathcal{HV}(a,b;1)}(\lambda,\alpha,\beta,\gamma):=\C[t]
$
by
\begin{eqnarray}
\label{action-hvab1-L-new} L_{i,m}\cdot f(t) &\!\!\!=\!\!\!& \sum_{s=0}^{+\infty}\binom{m}{s}\lambda^i\beta^{m-s}(t-i\alpha+(m-s)\alpha\beta^{-1})f^{(s)}(t-i),
\end{eqnarray}
\begin{eqnarray}
\label{action-hvab1-H-new} H_{i,m}\cdot f(t) &\!\!\!=\!\!\!&
\left\{
\begin{aligned}
&\sum_{s=0}^{+\infty}\binom{m}{s}\lambda^i \beta^{m-s}\gamma f^{(s)}(t-a-i), && \mbox{if\ } b=0,\\
&\sum_{s=0}^{+\infty}\sum_{\ell=0}^{m-s}\frac{(m-s)!}{(m-s-\ell)!}\binom{m}{s}\frac{a\lambda^i\beta^{m-s-\ell}\gamma}{(a+i)^{\ell+1}}f^{(s)}(t-a-i), && \mbox{if\ } b=1, a\ne0,\\
&0, && \mbox{otherwise},
\end{aligned}\right.
\end{eqnarray}
where $i\in\Z$, $m \in \Z_+$ and $f(t)\in\C[t]$.
One can easily check that $\Omega_{\mathcal{HV}(a,b;1)}(\lambda, \alpha, \beta, \gamma)$ is a $\mathcal{U}(\mathfrak h)$-module, which is free of rank one.
Furthermore, we find that it is in fact a $\mathcal{HV}(a,b;1)$-module.

\begin{prop}\label{hvab1-module}
The space $\Omega_{\mathcal{HV}(a,b;1)}(\lambda, \alpha, \beta, \gamma)$ is a $\mathcal{HV}(a,b;1)$-module under the actions \eqref{action-hvab1-L-new} and \eqref{action-hvab1-H-new}.
\end{prop}

\ni{\it Proof.}\ \
Let $V=\Omega_{\mathcal{HV}(a,b;1)}(\lambda, \alpha, \beta, \gamma)$.
Recall that $\{L_{i,0},L_{i,1},H_{i,0}\,|\, i\in\Z\}$ is a generating set of $\mathcal{HV}(a,b;1)$ (cf.~\eqref{gen-set}),
and both $\mathcal{W}(1)=\langle L_{i,0},L_{i,1}\,|\,i\in\Z\rangle$
and $\mathfrak{hv}(a,b)=\langle L_{i,0},H_{i,0}\,|\,i\in\Z\rangle$ are subalgebras of $\mathcal{HV}(a,b;1)$.
By Lemmas~\ref{We-module} and \ref{hvab-module}, it is enough to check the commutation relations between $L_{i,1}$ and $H_{j,0}$ on $V$.

We only check the case $b=1$, $a\ne0$; other cases are trivial or can be easily checked.
On the one hand, by \eqref{action-hvab1-H-new}, we have
\begin{eqnarray*}
\frac{1}{\lambda^{i+j}}[L_{i,1},H_{j,0}]\cdot f(t) &\!\!\!=\!\!\!& \frac{1}{\lambda^{i+j}}(a+i+j)H_{i+j,1}\cdot f(t) - \frac{1}{\lambda^{i+j}} H_{i+j,0}\cdot f(t) \\
&\!\!\!=\!\!\!& a\gamma\left(\left(\beta+\frac{1}{a+i+j}\right)f(t-a-i-j)+f'(t-a-i-j)\right) \\
&& -\frac{a\gamma}{a+i+j}f(t-a-i-j)\\
&\!\!\!=\!\!\!& a\gamma(\beta f(t-a-i-j)+f'(t-a-i-j)).
\end{eqnarray*}
On the other hand, by  \eqref{action-hvab1-L-new} and \eqref{action-hvab1-H-new}, we have
\begin{eqnarray*}
\frac{1}{\lambda^{i+j}}L_{i,1}\cdot H_{j,0}\cdot f(t) &\!\!\!=\!\!\!& \frac{1}{\lambda^{i}}L_{i,1}\cdot \left(\frac{a\gamma}{a+j}f(t-a-j)\right)\\
&\!\!\!=\!\!\!& \frac{a\gamma}{a+j}\left((\beta t-i\alpha\beta+\alpha)f(t-a-i-j)+(t-i\alpha)f'(t-a-i-j)\right),
\end{eqnarray*}
and
\begin{eqnarray*}
\frac{1}{\lambda^{i+j}}H_{j,0}\cdot L_{i,1}\cdot f(t) &\!\!\!=\!\!\!& \frac{1}{\lambda^{j}}H_{j,0}\cdot \left((\beta t-i\alpha\beta+\alpha)f(t-i)+(t-i\alpha)f'(t-i)\right)\\
&\!\!\!=\!\!\!& \frac{a\gamma}{a+j}\big((\beta (t-a-j)-i\alpha\beta+\alpha)f(t-a-i-j)\\
&& +(t-a-j-i\alpha)f'(t-a-i-j)\big).
\end{eqnarray*}
The above three formulas show that $[L_{i,1},H_{j,0}]\cdot f(t) = L_{i,1}\cdot H_{j,0}\cdot f(t) - H_{j,0}\cdot L_{i,1}\cdot f(t)$.
\QED

\begin{rema}\label{remark-case1}\rm
From \eqref{action-hvab1-L-new} and \eqref{action-hvab1-H-new}, we see clearly that
$\Omega_{\mathcal{HV}(a,b;1)}(\lambda,\alpha,\beta,\gamma)$ has four free parameters if $b=0$, or $b=1, a\ne0$,
and three free parameters if otherwise (as listed in Table~2).
\end{rema}

\subsection{\normalsize{Classification of free $\mathcal{U}(\mathfrak h)$-modules of rank one over $\mathcal{HV}(a,b;1)$}}

Now, we give the classification of free $\mathcal{U}(\mathfrak h)$-modules of rank one over $\mathcal{HV}(a,b;1)$.

\begin{theo}\label{thm-hvab1}
Let $\mathfrak{g}=\mathcal{HV}(a,b;1)$ and $\mathfrak h=\C L_{0,0}\in \mathfrak{g}$.
Any free $\mathcal{U}(\mathfrak h)$-module of rank one over $\mathfrak{g}$ is isomorphic to
$\Omega_{\mathfrak{g}}(\lambda,\alpha,\beta,\gamma)$ with actions \eqref{action-hvab1-L-new} and \eqref{action-hvab1-H-new} for some $\lambda\in\C^*$ and $\alpha, \beta, \gamma\in\C$.
\end{theo}

\ni{\it Proof.}\ \
Note first that the actions \eqref{action-hvab1-L-new} and \eqref{action-hvab1-H-new} are respectively equivalent to
\begin{eqnarray}
\label{action-hvab1-L} L_{i,m}\cdot t^k &\!\!\!=\!\!\!& \sum_{s=0}^{{\rm min}\{m,\,k\}}s!\binom{m}{s}\binom{k}{s}\lambda^i\beta^{m-s-1}((m-s)\alpha - i\alpha\beta + \beta t)(t-i)^{k-s},\\
\label{action-hvab1-H} H_{i,m}\cdot t^k &\!\!\!=\!\!\!&
\left\{
\begin{aligned}
&\sum_{s=0}^{{\rm min}\{m,k\}}s!\binom{m}{s}\binom{k}{s}\lambda^i \beta^{m-s}\gamma(t-a-i)^{k-s}, && \mbox{if\ } b=0,\\
&\sum_{s=0}^{{\rm min}\{m,k\}}\sum_{\ell=0}^{m\!-\!s}\frac{s!(m-s)!}{(m\!-\!s\!-\!\ell)!}\binom{m}{s}\!\!\binom{k}{s}\frac{a\lambda^i\beta^{m-s-\ell}\gamma}{(a+i)^{\ell+1}}(t\!-\!a\!-\!i)^{k-s}, && \mbox{if\ } b=1, a\ne0,\\
&0, && \mbox{otherwise}.
\end{aligned}\right.
\end{eqnarray}
Let $M$ be a free $\mathcal{U}(\mathfrak h)$-module of rank one over $\mathfrak{g}$.
As we mentioned in the Introduction, both $\mathcal{W}(\epsilon)$ and $\mathfrak{hv}(a,b)$ are subalgebras of $\mathfrak{g}$,
and they share the same abelian subalgebra $\mathfrak h$.
By viewing $M$ as a module over both $\mathcal{W}(\epsilon)$ and $\mathfrak{hv}(a,b)$, from Lemmas~\ref{We-module} and \ref{hvab-module},
we may assume that $M=\C[t]$ and there exist some $\lambda\in\C^*$ and $\alpha,\beta,\gamma\in\C$ such that the action of $L_{i,m}$ on $t^k$ has form \eqref{action-hvab1-L},
and
\begin{equation}\label{Hi0t^k}
H_{i,0}\cdot t^k = \phi_i(a,b)\lambda^i\gamma (t-a-i)^k,
\end{equation}
where $\phi_i(a,b)$ is given by \eqref{func-phiiab}.
The formula \eqref{action-hvab1-H} will follow from the following Lemmas~\ref{hvab1-action-on-t^k} and \ref{hvab1-action-on-1},
and the relation $\ell !\binom{m-s}{\ell}=\frac{(m-s)!}{(m-s-\ell)!}$.
\QED

\begin{lemm}\label{hvab1-action-on-t^k}
The action of $H_{i,m}$ on $t^k$ is a combination of the actions of $H_{i,j}$ ($m-{\rm min}\{m,k\} \leq j \leq m$) on $1$,
and more precisely we have
\begin{equation*}\label{claim-hvab1-relation}
H_{i,m}\cdot t^k = \sum_{s=0}^{{\rm min}\{m,k\}}s!\binom{m}{s}\binom{k}{s}(t-a-i)^{k-s}H_{i,m-s}\cdot 1.
\end{equation*}
\end{lemm}
\ni{\it Proof.}\ \
The proof is similar to that of \cite[Lemma~3.3]{XMDZ2023} and is omitted.
The main strategy is to apply induction on $k$, together with some combinatorial techniques, including Pascal's triangle.
\QED

\begin{lemm}\label{hvab1-action-on-1}
The action of $H_{i,m}$ on $1$ is as follows:
\begin{eqnarray*}
H_{i,m}\cdot 1 &\!\!\!=\!\!\!&
\left\{
\begin{aligned}
&\lambda^i\beta^m\gamma, && \mbox{if\ } b=0,\\
&a\lambda^i\gamma\sum_{\ell=0}^{m}\ell!\binom{m}{\ell}\frac{\beta^{m-\ell}}{(a+i)^{\ell+1}}, && \mbox{if\ } b=1, a\ne0,\\
&0, && \mbox{otherwise}.
\end{aligned}\right.
\end{eqnarray*}
\end{lemm}

\ni{\it Proof.}\ \
Denote $G_{i,m}(t)=H_{i,m}\cdot 1$. %Then \eqref{Hi0t^k} with $k=0$ gives $G_{i,0}(t)=\phi_i(a,b)\lambda^i\gamma$.
We shall determine $G_{i,m}(t)$ according to the classification \eqref{func-phiiab} of $\phi_i(a,b)$.
The formulas \eqref{action-hvab1-L} with $k=0$ and \eqref{Hi0t^k} with $k=0,1$ will be used frequently:
\begin{equation}\label{L-k=0}
L_{i,m}\cdot 1 = \lambda^i\beta^{m-1}(m\alpha - i\alpha\beta + \beta t),
\end{equation}
\begin{equation}\label{H-k=0,1}
H_{i,0}\cdot 1 = \phi_i(a,b)\lambda^i\gamma, \quad H_{i,0}\cdot t = \phi_i(a,b)\lambda^i\gamma(t-a-i).
\end{equation}

{\bf Case 1:} $b=0$.

In this case, $\phi_i(a,b)=1$.
Applying $[L_{0,m},H_{i,0}]= (a+i)H_{i,m}$ on $1$, by \eqref{L-k=0} and \eqref{H-k=0,1},
we have %one can derive that
%$$
%(a+i)G_{i,m}(t)=(a+i)\lambda^i\beta^m\gamma.
%$$
\begin{eqnarray*}
(a+i)G_{i,m}(t) &\!\!\!=\!\!\!& L_{0,m}\cdot H_{i,0}\cdot 1 - H_{i,0}\cdot L_{0,m}\cdot 1\\
&\!\!\!=\!\!\!&  L_{0,m}\cdot (\lambda^i\gamma) - H_{i,0}\cdot (\beta^mt + m\alpha\beta^{m-1})\\
&\!\!\!=\!\!\!&  \lambda^i\gamma(\beta^mt + m\alpha\beta^{m-1}) - \beta^m(\lambda^i\gamma(t-a-i)) - m\alpha\beta^{m-1}(\lambda^i\gamma)\\
&\!\!\!=\!\!\!& (a+i)\lambda^i\beta^m\gamma.
\end{eqnarray*}
Recall that $0 \leq Re(a)< 1$ (Convention~\ref{conv-iso}).
If $a\ne0$, then $a+i\ne0$ for $i\in\Z$. Then the above equation implies that $G_{i,m}(t) = \lambda^i\beta^m\gamma$, as desired.
If $a=0$, then the above equation implies that $G_{i,m}(t) = \lambda^i\beta^m\gamma$ for $i\ne0$.
Furthermore, applying $[L_{-1,m}, H_{1,0}]=H_{0,m}$ on $1$, one can derive that $G_{0,m}(t) = \beta^m\gamma$, as desired.

{\bf Case 2:} $b=1, a \neq 0$.

In this case, we have $\phi_i(a,b)=\frac{a}{a+i}$. We need to prove
\begin{equation}\label{case2}
G_{i,m}(t)=a\lambda^i\gamma\sum_{\ell=0}^{m}\ell!\binom{m}{\ell}\frac{\beta^{m-\ell}}{(a+i)^{\ell+1}}.
\end{equation}
Our proof is by induction on $m$. The case for $m=0$ is given by the first formula in \eqref{H-k=0,1}.
Let $m\ge 1$. Applying $[L_{0,m},H_{i,0}] = (a+i)H_{i,m} - m H_{i,m-1}$ on $1$, by \eqref{L-k=0} and \eqref{H-k=0,1}, one can derive that
\begin{equation}\label{case2-relation}
(a+i)G_{i,m}(t) - mG_{i,m-1}(t)=a\lambda^i\beta^m\gamma. 
\end{equation}
By \eqref{case2-relation} and induction hypothesis, we have
\begin{eqnarray*}
G_{i,m}(t) &\!\!\!=\!\!\!& \frac{a}{a+i}\lambda^i\beta^m\gamma + \frac{m}{a+i}G_{i,m-1}(t)\\
 &\!\!\!=\!\!\!& \frac{a}{a+i}\lambda^i\beta^m\gamma + \frac{ma}{a+i}\lambda^i\gamma\sum_{\ell=0}^{m-1}\ell!\binom{m-1}{\ell}\frac{\beta^{m-1-\ell}}{(a+i)^{\ell+1}}\\
 &\!\!\!=\!\!\!& a\lambda^i\gamma\frac{\beta^m}{a+i} + a\lambda^i\gamma\sum_{\ell=0}^{m-1}(\ell+1)!\binom{m}{\ell+1}\frac{\beta^{m-1-\ell}}{(a+i)^{\ell+2}}\\
 &\!\!\!=\!\!\!& a\lambda^i\gamma\frac{\beta^m}{a+i} + a\lambda^i\gamma\sum_{\ell=1}^{m}\ell!\binom{m}{\ell}\frac{\beta^{m-\ell}}{(a+i)^{\ell+1}}\\
 &\!\!\!=\!\!\!& a\lambda^i\gamma\sum_{\ell=0}^{m}\ell!\binom{m}{\ell}\frac{\beta^{m-\ell}}{(a+i)^{\ell+1}}.
\end{eqnarray*}
Hence, \eqref{case2}  holds.

{\bf Case 3:} $b=1, a = 0$.

In this case, we have $\phi_i(a,b)=\delta_{i,0}$.
Applying $[L_{0,m},H_{i,0}] = iH_{i,m} - mH_{i,m-1}$ on $1$, by \eqref{L-k=0} and \eqref{H-k=0,1}, we obtain
\begin{equation}\label{case3-relation}
iG_{i,m}(t) - mG_{i,m-1}(t)=0.
\end{equation}
Taking $i=0$ and $m=1$ in \eqref{case3-relation}, we see that $\gamma=0$.
By \eqref{Hi0t^k}, we have $H_{i,0}\cdot t^k = 0$ for $i\in\Z$.
In particular, $G_{i,0}(t)=H_{i,0}\cdot 1=0$.
Then, by \eqref{case3-relation}, one can inductively see that $G_{i,m}(t)=0$ for $i\ne0$.
Taking $i=0$ in \eqref{case3-relation} again, we have $mG_{0,m-1}(t)=0$. This implies that $G_{0,m}(t)=0$ for $m\in\Z_+$, as desired.

{\bf Case 4:} $b \neq 0,1$.

In this case, we have $\phi_i(a,b)=0$. The first formula in \eqref{H-k=0,1} gives $G_{i,0}(t)=H_{i,0}\cdot 1=0$.
Applying $[L_{0,m},H_{i,0}] = (a+i)H_{i,m} - bmH_{i,m-1}$ on $1$, we obtain
\begin{equation}\label{case4-relation}
(a+i)G_{i,m}(t) - bmG_{i,m-1}(t)=0.
\end{equation}
Recall that $0 \leq Re(a)< 1$ (Convention~\ref{conv-iso}).
If $a\ne0$, then $a+i\ne0$ for $i\in\Z$. The formula \eqref{case4-relation} inductively implies that $G_{i,m}(t)=0$.
If $a=0$, one can derive $G_{i,m}(t)=0$ as in Case~3.
\QED

\subsection{\normalsize{Simplicity and isomorphism classification of $\mathcal{HV}(a,b;1)$-modules}}

Combining Lemmas~\ref{we-sim&iso} and \ref{hvab-sim&iso}, we immediately obtain the following results.

\begin{theo}\label{hvab1-sim&iso}
Let $\mathfrak{g}=\mathcal{HV}(a,b;1)$. The following statements hold:
\begin{itemize}\parskip-3pt
\item[{\rm(1)}] $\Omega_{\mathfrak{g}}(\lambda,\alpha,\beta,\gamma)$ is simple if and only if $\alpha \neq 0\ {\rm or}\ \gamma\neq 0$ (if $\gamma$ exists);
\item[{\rm(2)}] $\Omega_{\mathfrak{g}}(\lambda,0,\beta,0)$ has a unique proper submodule $t\Omega_{\mathfrak{g}}(\lambda,0,\beta,0)\cong\Omega_{\mathfrak{g}}(\lambda,1,\beta,0)$, and the quotient $\Omega_{\mathfrak{g}}(\lambda,0,\beta,0)/t\Omega_{\mathfrak{g}}(\lambda,0,\beta,0)$ is a one-dimensional trivial $\mathfrak{g}$-module;
\item[{\rm(3)}] $\Omega_{\mathfrak{g}}(\lambda_1,\alpha_1,\beta_1,\gamma_1)\cong\Omega_{\mathfrak{g}}(\lambda_2,\alpha_2,\beta_2,\gamma_2)$ if and only if they have the same parameters.
\end{itemize}
\end{theo}

\section{\large{Free $\mathcal{U}(\mathfrak h)$-modules of rank one over $\mathcal{HV}(a,b;-1)$}}

Recall again that $\mathfrak h=\C L_{0,0}\in \mathcal{HV}(a,b;-1)$.
In this section, we completely classify the free $\mathcal{U}(\mathfrak h)$-modules of rank one over $\mathcal{HV}(a,b;-1)$.
One will see that the module structures of $\mathcal{HV}(a,b;-1)$ are more varied and complex than those of $\mathcal{W}(-1)$.
In particular, if $b=1$, there will appear unusual modules which admit infinitely many free parameters 
(Remark~\ref{remark-case-1} and Theorem~\ref{thm-hvab-1}),
and thus give an affirmative answer to the question we proposed in Subsection~2.1.

\subsection{\normalsize{Construction of free $\mathcal{U}(\mathfrak h)$-modules of rank one over $\mathcal{HV}(a,b;-1)$}}

Let $\lambda\in\C^*$, $\alpha,\beta,\gamma\in\C$, and $\bm{\kappa}=\{(\kappa_i)_{i\in\Z} \mid \kappa_i\in\C\}$.
Define the action of $\mathcal{HV}(a,b;-1)$ on the vector space of polynomials in one variable
$
\Omega_{\mathcal{HV}(a,b;-1)}(\lambda,\alpha,\beta,\gamma,\bm{\kappa}):=\C[t]
$
by
\begin{eqnarray}
\label{action-hvab-1-L-new} L_{i,m}\cdot f(t) &\!\!\!=\!\!\!& \sum_{s=0}^{+\infty}\binom{-m}{s}\lambda^i\beta^{m+s}(t-i\alpha-(m+s)\alpha\beta)f^{(s)}(t-i),\\
\label{action-hvab-1-H-new} H_{i,m}\cdot f(t) &\!\!\!=\!\!\!&
\left\{
\begin{aligned}
&\sum_{s=0}^{+\infty}\binom{-m}{s}\lambda^i\beta^{m+s}\gamma f^{(s)}(t-a-i), && \mbox{if\ } b=0,\\
&\sum_{s=0}^{+\infty}\binom{-m}{s}\varphi_{i,m+s}^{(a)}(\lambda,\beta,\gamma,\kappa_i)f^{(s)}(t-a-i),&& \mbox{if\ } b=1, a\ne0,\\
&\sum_{s=0}^{+\infty}\binom{-m}{s}\psi_{i,m+s}(\gamma,\kappa_i)f^{(s)}(t-i), && \mbox{if\ } b = 1, a = 0,\\
&0, && \mbox{otherwise},
\end{aligned}\right.
\end{eqnarray}
where $i\in\Z$, $m \in \Z_+$, $f(t)\in\C[t]$, and $\varphi_{i,m}^{(a)}(\lambda,\beta,\gamma,\kappa_i)$ and $\psi_{i,m}(\gamma,\kappa_i)$ are respectively defined by
\begin{eqnarray}
\nonumber \varphi_{i,m}^{(a)}(\lambda,\beta,\gamma,\kappa_i) &\!\!\!=\!\!\!& \delta_{m,0}\frac{a}{a+i}\lambda^i\gamma + \frac{(-1)^{m-1}}{(m-1)!}(a+i)^{m-1}\kappa_i  \\
\label{fun-varphi} &\!\!\!\!\!\!&+\sum_{\ell=0}^{m-2}\frac{(m-\ell-2)!}{(m-1)!}(-1)^{\ell}(a+i)^{\ell}a\lambda^i\beta^{m-\ell-1}\gamma, \\
\label{fun-psi} \psi_{i,m}(\gamma,\kappa_i) &\!\!\!=\!\!\!& \delta_{i,0}\delta_{m,0}\gamma + \frac{(-i)^{m-1}}{(m-1)!}\kappa_i.
\end{eqnarray}
One can easily check that $\Omega_{\mathcal{HV}(a,b;-1)}(\lambda,\alpha,\beta,\gamma,\bm{\kappa})$ is a $\mathcal{U}(\mathfrak h)$-module, which is free of rank one.
Similar to Proposition~\ref{hvab1-module}, one can further prove that it is a $\mathcal{HV}(a,b;-1)$-module.

\begin{prop}\label{hvab-1-module}
The space $\Omega_{\mathcal{HV}(a,b;-1)}(\lambda,\alpha,\beta,\gamma,\bm{\kappa})$ is a $\mathcal{HV}(a,b;-1)$-module under the actions
\eqref{action-hvab-1-L-new} and \eqref{action-hvab-1-H-new}.
\end{prop}

\ni{\it Proof.}\ \
Let $V=\Omega_{\mathcal{HV}(a,b;-1)}(\lambda,\alpha,\beta,\gamma,\bm{\kappa})$.
Similar to Proposition~\ref{hvab1-module}, we only check the most complicated case $b=1$, $a\ne0$.
Recall that $\{L_{i,0},L_{i,1},H_{i,0},H_{i,1}\,|\,i\in\Z\}$ is a generating set of $\mathcal{HV}(a,1;-1)$ (cf.~\eqref{gen-set}).
By Lemmas~\ref{We-module} and \ref{hvab-module}, we need to further check the following brackets on $V$:
$$
[L_{i,0},H_{j,1}],\ [L_{i,1},\ H_{j,0}],\ [L_{i,1},H_{j,1}],\ [H_{i,0},H_{j,1}],\ [H_{i,1},H_{j,1}].
$$
To save space, here, we only check the bracket $[L_{i,1},H_{j,0}]$ on $V$ (we leave the details of other brackets to interested readers).

On the one hand, by \eqref{action-hvab-1-H-new}, we have
\begin{eqnarray*}
[L_{i,1}, H_{j,0}]\cdot f(t)&\!\!\!=\!\!\!& (a+i+j)H_{i+j,1}\cdot f(t) + H_{i+j,2}\cdot f(t)=W(t)+X(t)+Y(t)+Z(t),
\end{eqnarray*}
where
\begin{eqnarray*}
W(t) &\!\!\!=\!\!\!& (a+i+j)\left(\sum_{s=0}^{+\infty}\frac{1}{s!}(a+i+j)^{s}\kappa_{i+j}f^{(s)}(t-a-i-j)\right),\\
X(t) &\!\!\!=\!\!\!& (a+i+j)\left(\sum_{s=0}^{+\infty}\sum_{\ell=0}^{s-1}\frac{(s-\ell-1)!}{s!}(-1)^{s+\ell}(a+i+j)^{\ell}a\lambda^{i+j}\beta^{s-\ell}\gamma f^{(s)}(t-a-i-j)\right),\\
Y(t) &\!\!\!=\!\!\!& \sum_{s=0}^{+\infty}(-1)^{s}(s+1)\left(\frac{(-1)^{s+1}}{(s+1)!}(a+i+j)^{s+1}\kappa_{i+j} \right)f^{(s)}(t-a-i-j),\\
Z(t) &\!\!\!=\!\!\!& \sum_{s=0}^{+\infty}(-1)^{s}(s+1)\left(\sum_{\ell=0}^{s}\frac{(s-\ell)!}{(s+1)!}(-1)^{\ell}(a+i+j)^{\ell}a\lambda^{i+j}\beta^{s-\ell+1}\gamma\right)f^{(s)}(t-a-i-j).
\end{eqnarray*}
Clearly, $W(t)+Y(t)=0$. In addition, by the following observation
$$
\sum_{\ell=0}^{s-1}\frac{(s-\ell-1)!}{s!}(-1)^{s+\ell}(a+i+j)^{\ell+1}\beta^{s-\ell}
+\sum_{\ell=0}^{s}\frac{(s-\ell)!}{s!}(-1)^{s+\ell}(a+i+j)^{\ell}\beta^{s-\ell+1}=(-1)^s \beta^{s+1},
$$
we see that
$$
X(t)+Z(t)=a\lambda^{i+j}\gamma\sum_{s=0}^{+\infty}(-1)^s \beta^{s+1} f^{(s)}(t-a-i-j).
$$
On the other hand, by \eqref{action-hvab-1-L-new} and \eqref{action-hvab-1-H-new}, we have
\begin{eqnarray*}
L_{i,1}\cdot H_{j,0}\cdot f(t) &\!\!\!=\!\!\!& L_{i,1}\cdot \left(\frac{a\lambda^j\gamma}{a+j}f(t-a-j)\right)\\
&\!\!\!=\!\!\!& \frac{a\lambda^{i+j}\gamma}{a+j}\sum_{s=0}^{+\infty}(-1)^s \beta^{s+1} (t-i\alpha-(1+s)\alpha\beta) f^{(s)}(t-a-i-j),
\end{eqnarray*}
and
\begin{eqnarray*}
H_{j,0}\cdot L_{i,1}\cdot f(t) &\!\!\!=\!\!\!& H_{j,0}\cdot \sum_{s=0}^{+\infty}(-1)^s \lambda^i \beta^{s+1} (t-i\alpha-(1+s)\alpha\beta) f^{(s)}(t-i)\\
&\!\!\!=\!\!\!& \frac{a\lambda^{i+j}\gamma}{a+j}\sum_{s=0}^{+\infty}(-1)^s \beta^{s+1} (t-a-j-i\alpha-(1+s)\alpha\beta) f^{(s)}(t-a-i-j),
\end{eqnarray*}
and thus
$$
L_{i,1}\cdot H_{j,0}\cdot f(t) - H_{j,0}\cdot L_{i,1}\cdot f(t)=a\lambda^{i+j}\gamma\sum_{s=0}^{+\infty}(-1)^s \beta^{s+1} f^{(s)}(t-a-i-j).
$$
Hence, we have proved that $[L_{i,1}, H_{j,0}]\cdot f(t)=L_{i,1}\cdot H_{j,0}\cdot f(t) - H_{j,0}\cdot L_{i,1}\cdot f(t)$.
\QED

\begin{rema}\label{remark-case-1}\rm
From \eqref{action-hvab-1-L-new}--\eqref{fun-psi}, we see that $\Omega_{\mathcal{HV}(a,b;-1)}(\lambda,\alpha,\beta,\gamma,\bm{\kappa})$ has four free parameters if $b=0$,
infinitely many free parameters if $b=1$, and three free parameters if otherwise (as listed in Table~2).
Hence, from the view point of the number of free parameters, the free $\mathcal{U}(\mathfrak h)$-modules over $\mathcal{HV}(a,1;-1)$ are distinctive.
\end{rema}

\subsection{\normalsize{Classification of free $\mathcal{U}(\mathfrak h)$-modules of rank one over $\mathcal{HV}(a,b;-1)$}}

Now, we classify free $\mathcal{U}(\mathfrak h)$-modules of rank one over $\mathcal{HV}(a,b;-1)$.

\begin{theo}\label{thm-hvab-1}
Let $\mathfrak{g}=\mathcal{HV}(a,b;-1)$ and $\mathfrak h=\C L_{0,0}\in \mathfrak{g}$.
Any free $\mathcal{U}(\mathfrak h)$-module of rank one over $\mathfrak{g}$ is isomorphic to
$\Omega_{\mathfrak{g}}(\lambda,\alpha,\beta,\gamma,\bm{\kappa})$ with actions \eqref{action-hvab-1-L-new} and \eqref{action-hvab-1-H-new}
for some $\lambda\in\C^*$, $\alpha, \beta, \gamma\in\C$, and $\bm{\kappa}=\{(\kappa_i)_{i\in\Z} \mid \kappa_i\in\C\}$.
\end{theo}

\ni{\it Proof.}\ \
For $m\in\Z$ and $n\in\Z_+$, recall that
$
\binom{-m}{n}=(-1)^n\binom{m+n-1}{n} \mbox{\ if\ } m>0.
$
We note that the actions \eqref{action-hvab-1-L-new} and \eqref{action-hvab-1-H-new} are respectively equivalent to
\begin{eqnarray}
\label{action-hvab-1-L} L_{i,m}\cdot t^k &\!\!\!=\!\!\!& \sum_{s=0}^{k}(-1)^ss!\binom{m+s-1}{s}\binom{k}{s}\lambda^i\beta^{m+s}(t-i\alpha-(m+s)\alpha\beta)(t-i)^{k-s},\\
\label{action-hvab-1-H}H_{i,m}\cdot t^k &\!\!\!=\!\!\!&
\left\{
\begin{aligned}
&\sum_{s=0}^{k}(-1)^ss!\binom{m+s-1}{s}\binom{k}{s}\lambda^i\beta^{m+s}\gamma(t-a-i)^{k-s}, && \mbox{if\ } b=0,\\
&\sum_{s=0}^{k}(-1)^ss!\binom{m\!+\!s\!-\!1}{s}\!\!\binom{k}{s}\varphi_{i,m+s}^{(a)}(\lambda,\!\beta,\!\gamma,\!\kappa_i)(t\!-\!a\!-\!i)^{k-s},&& \mbox{if\ } b=1, a\ne0,\\
&\sum_{s=0}^{k}(-1)^ss!\binom{m\!+\!s\!-\!1}{s}\!\!\binom{k}{s}\psi_{i,m+s}(\gamma,\kappa_i)(t-i)^{k-s}, && \mbox{if\ } b = 1, a = 0,\\
&0, && \mbox{otherwise}.
\end{aligned}\right.
\end{eqnarray}
Let $M$ be a free $\mathcal{U}(\mathfrak h)$-module of rank one over $\mathfrak{g}$.
Similar to Theorem~\ref{thm-hvab1},
by Lemmas~\ref{We-module} and \ref{hvab-module},
we may assume that $M=\C[t]$ and there exist some $\lambda\in\C^*$ and $\alpha,\beta,\gamma\in\C$ such that the action of $L_{i,m}$ on $t^k$ has form \eqref{action-hvab-1-L},
and
\begin{equation}\label{Hi0t^k-case-1}
H_{i,0}\cdot t^k = \phi_i(a,b)\lambda^i\gamma (t-a-i)^k,
\end{equation}
where $\phi_i(a,b)$ is given by \eqref{func-phiiab}.
The formula \eqref{action-hvab-1-H} will follow from the following Lemmas~\ref{hvab-1-action-on-t^k} and \ref{hvab-1-action-on-1}
(note that the parameter set $\bm{\kappa}$ will appear in Lemma~\ref{hvab-1-action-on-1}).
\QED

\begin{lemm}\label{hvab-1-action-on-t^k}
The action of $H_{i,m}$ on $t^k$ is a combination of the actions of $H_{i,j}$ ($m \leq j \leq m+k$) on $1$,
and more precisely we have
\begin{equation*}\label{claim-hvab-1-relation}
H_{i,m}\cdot t^k = \sum_{s=0}^{k}(-1)^ss!\binom{m+s-1}{s}\binom{k}{s}(t-a-i)^{k-s}H_{i,m+s}\cdot 1.
\end{equation*}
\end{lemm}
\ni{\it Proof.}\ \
The proof is similar to that of \cite[Lemma~3.3]{XMDZ2023} and is omitted.
\QED

\begin{lemm}\label{hvab-1-action-on-1}
The action of $H_{i,m}$ on $1$ is as follows:
\begin{eqnarray*}
H_{i,m}\cdot 1 &\!\!\!=\!\!\!&
\left\{
\begin{aligned}
&\lambda^i\beta^m\gamma, && \mbox{if\ } b=0,\\
&\varphi_{i,m}^{(a)}(\lambda,\beta,\gamma,\kappa_i), && \mbox{if\ } b=1, a\ne0,\\
&\psi_{i,m}(\gamma,\kappa_i), && \mbox{if\ } b=1, a=0,\\
&0, && \mbox{otherwise},
\end{aligned}\right.
\end{eqnarray*}
where $\varphi_{i,m}^{(a)}(\lambda,\beta,\gamma,\kappa_i)$ and $\psi_{i,m}(\gamma,\kappa_i)$ are defined by
\eqref{fun-varphi} and \eqref{fun-psi}, respectively.
\end{lemm}

\ni{\it Proof.}\ \
Denote $G_{i,m}(t)=H_{i,m}\cdot 1$.
Similar to Lemma~\ref{hvab1-action-on-1}, we determine $G_{i,m}(t)$ according to the classification \eqref{func-phiiab} of $\phi_i(a,b)$.
The following formulas (cf.~\eqref{action-hvab-1-L} and \eqref{Hi0t^k-case-1}) will be used frequently:
\begin{equation}\label{L-k=0-case-1}
L_{i,m}\cdot 1 = \lambda^i\beta^{m}(t-i\alpha-m\alpha\beta),
\end{equation}
\begin{equation}\label{H-k=0,1-case-1}
H_{i,0}\cdot 1 = \phi_i(a,b)\lambda^i\gamma, \quad H_{i,0}\cdot t = \phi_i(a,b)\lambda^i\gamma(t-a-i).
\end{equation}

{\bf Case 1:} $b=0$.

In this case, $\phi_i(a,b)=1$. The proof is similar to that of Case~1 in Lemma~\ref{hvab1-action-on-1} and is omitted.

{\bf Case 2:} $b=1, a \neq 0$.

In this case, we have  $\phi_i(a,b)=\frac{a}{a+i}$.
The first formula in \eqref{H-k=0,1-case-1} gives $H_{i,0}\cdot 1=\varphi_{i,0}^{(a)}(\lambda,\beta,\gamma,\kappa_i)$ (which is in fact not depend on $\beta$ and $\kappa_i$).
Thus, we only need to prove that, for $m\ge1$, there exists $\bm{\kappa}=\{(\kappa_i)_{i\in\Z} \mid \kappa_i\in\C\}$ such that
\begin{equation}\label{cliam-4-case2-1}
G_{i,m}(t) = \frac{(-1)^{m-1}}{(m-1)!}(a+i)^{m-1}\kappa_i + \sum_{\ell=0}^{m-2}\frac{(m-\ell-2)!}{(m-1)!}(-1)^{\ell}(a+i)^{\ell}a\lambda^i\beta^{m-\ell-1}\gamma.
\end{equation}

Let us first show that $G_{i,m}(t)\in\C$. Applying relations
\begin{eqnarray*}
\nonumber [L_{0,m}, H_{i,0}] &\!\!\!=\!\!\!& (a+i)H_{i,m} + m H_{i,m+1},\\
\nonumber [L_{-i,0}, H_{i,m}] &\!\!\!=\!\!\!& aH_{0,m} + mH_{0,m+1}, \\
\nonumber [L_{i,0}, H_{0,m}] &\!\!\!=\!\!\!& (a+i)H_{i,m} + mH_{i,m+1}
\end{eqnarray*}
respectively on $1$, by \eqref{L-k=0-case-1}, \eqref{H-k=0,1-case-1} and Lemma~\ref{hvab-1-action-on-t^k}, one can derive that
\begin{eqnarray}
\label{action-case2-Gimt} (a+i)G_{i,m}(t)+mG_{i,m+1}(t) &\!\!\!=\!\!\!& a\lambda^i\beta^m\gamma, \\
\nonumber aG_{0,m}(t) + mG_{0,m+1}(t) &\!\!\!=\!\!\!& \lambda^{-i}\big((t+i\alpha)G_{i,m}(t+i)-(t-a-i+i\alpha)G_{i,m}(t) \\
\label{action-case2-Gimt-1} &\!\!\!\!\!\!& + mG_{i,m+1}(t)\big),\\
\nonumber (a+i)G_{i,m}(t) + mG_{i,m+1}(t) &\!\!\!=\!\!\!& \lambda^{i}\big((t-i\alpha)G_{0,m}(t-i)-(t-a-i\alpha)G_{0,m}(t) \\
\label{action-case2-Gimt-2} &\!\!\!\!\!\!& + mG_{0,m+1}(t)\big).
\end{eqnarray}
In particular, taking $i=0$ in \eqref{action-case2-Gimt}, we have
\begin{eqnarray}\label{action-case2-Gimt-i=0}
aG_{0,m}(t)+mG_{0,m+1}(t) = a\beta^m\gamma.
\end{eqnarray}
By \eqref{action-case2-Gimt}, \eqref{action-case2-Gimt-1} and \eqref{action-case2-Gimt-i=0}, we see that
$\lambda^{-i}(t+i\alpha)(G_{i,m}(t+i) - G_{i,m}(t))=0$, which implies that $G_{i,m}(t+i)= G_{i,m}(t)$,
and thus $G_{i,m}(t)\in\C$ for $i\ne 0$.
Similarly, by \eqref{action-case2-Gimt}, \eqref{action-case2-Gimt-2} and \eqref{action-case2-Gimt-i=0}, we see that
$\lambda^{i}(t-i\alpha)(G_{0,m}(t-i)-G_{0,m}(t))=0$, which implies that $G_{0,m}(t-i) = G_{0,m}(t)$, and thus $G_{0,m}(t)\in\C$, as desired.

Now, let $\bm{\kappa}=\{(\kappa_i)_{i\in\Z} \mid \kappa_i=G_{i,1}(t)\in\C\}$. Then, we see that \eqref{cliam-4-case2-1} holds for $m=1$.
Assume that \eqref{cliam-4-case2-1} holds for $m$ ($\ge1$). By \eqref{action-case2-Gimt}, we have
\begin{eqnarray*}
 &\!\!\!\!\!\!& G_{i,m+1}(t) \\
 &\!\!\!=\!\!\!& \frac{a}{m}\lambda^i\beta^m\gamma - \frac{(a+i)}{m}G_{i,m}(t)\\
&\!\!\!=\!\!\!& \frac{a}{m}\lambda^i\beta^m\gamma - \frac{(a+i)}{m}\left(\frac{(-1)^{m-1}}{(m-1)!}(a+i)^{m-1}\kappa_i + \sum_{\ell=0}^{m-2}\frac{(m-\ell-2)!}{(m-1)!}(-1)^{\ell}(a+i)^{\ell}a\lambda^i\beta^{m-\ell-1}\gamma\right)\\
&\!\!\!=\!\!\!& \frac{a}{m}\lambda^i\beta^m\gamma +\frac{(-1)^{m}}{m!}(a+i)^{m}\kappa_i + \sum_{\ell=0}^{m-2}\frac{(m-\ell-2)!}{m!}(-1)^{\ell+1}(a+i)^{\ell+1}a\lambda^i\beta^{m-\ell-1}\gamma\\
&\!\!\!=\!\!\!& \frac{a}{m}\lambda^i\beta^m\gamma +\frac{(-1)^{m}}{m!}(a+i)^{m}\kappa_i + \sum_{\ell=1}^{m-1}\frac{(m-\ell-1)!}{m!}(-1)^{\ell}(a+i)^{\ell}a\lambda^i\beta^{m-\ell}\gamma\\
&\!\!\!=\!\!\!&\frac{(-1)^{m}}{m!}(a+i)^{m}\kappa_i + \sum_{\ell=0}^{m-1}\frac{(m-\ell-1)!}{m!}(-1)^{\ell}(a+i)^{\ell}a\lambda^i\beta^{m-\ell}\gamma.
\end{eqnarray*}
Hence, \eqref{cliam-4-case2-1} holds for all $m\ge1$.

{\bf Case 3:} $b=1, a = 0$.

In this case, we have $\phi_i(a,b)=\delta_{i,0}$.
The first formula in \eqref{H-k=0,1-case-1} gives $H_{i,0}\cdot 1=\delta_{i,0} \gamma=\psi_{i,0}(\gamma,\kappa_i)$
(the second equality is given by \eqref{fun-psi} with $m=1$).
Thus, we only need to prove that, for $m\ge1$, there exists $\bm{\kappa}=\{(\kappa_i)_{i\in\Z} \mid \kappa_i\in\C\}$ such that
\begin{equation}\label{cliam-4-case3-1}
G_{i,m}(t) = \frac{(-i)^{m-1}}{(m-1)!}\kappa_i.
\end{equation}
Using the same arguments as in Case~2, one can first show that $G_{i,m}(t)\in\C$.
Let $\bm{\kappa}=\{(\kappa_i)_{i\in\Z} \mid \kappa_i=G_{i,1}(t)\in\C\}$. We see that \eqref{cliam-4-case3-1} holds for $m=1$.
By relation $iG_{i,m}(t)+mG_{i,m+1}(t)=0$ (cf.~\eqref{action-case2-Gimt} with $a=0$), one can easily prove \eqref{cliam-4-case3-1} by induction on $m$.

{\bf Case 4:} $b \neq 0,1$.

In this case, we have $\phi_i(a,b)=0$.
By \eqref{Hi0t^k-case-1}, we see that the action of $H_{i,0}$ on $M$ is trivial.
Applying $[L_{j,m},H_{i,0}] = (a+i+bj)H_{i+j,m} + bm H_{i+j,m+1}$ on $1$, we obtain
\begin{equation}\label{LjmHi0}
(a+i+bj)G_{i+j,m}(t)+bm G_{i+j,m+1}(t)=0.
\end{equation}
Making the replacement $i\rightarrow i+1$, $j\rightarrow j-1$ in \eqref{LjmHi0}, we have
\begin{equation}\label{LjmHi0-new}
(a+i+1+b(j-1))G_{i+j,m}(t)+bm G_{i+j,m+1}(t)=0.
\end{equation}
By \eqref{LjmHi0} and \eqref{LjmHi0-new}, one can easily see that $(1-b)G_{i+j,m}(t)=0$.
Since $b\ne 1$, we must have $G_{i,m}(t)=0$, as desired.
\QED

\subsection{\normalsize{Simplicity and isomorphism classification of $\mathcal{HV}(a,b;-1)$-modules}}

Finally, we give the simplicity and isomorphism classification of the $\mathcal{HV}(a,b;-1)$-modules
$\Omega_{\mathcal{HV}(a,b;-1)}(\lambda,\alpha,\beta,\gamma,\bm{\kappa})$ constructed in Subsection~4.1.

\begin{theo}\label{hvab-1-sim&iso}
Let $\mathfrak{g}=\mathcal{HV}(a,b;-1)$ and $\bm{0}=\{(\kappa_i)_{i\in\Z} \mid \kappa_i=0\}$.
 The following statements hold:
\begin{itemize}\parskip-3pt
\item[{\rm(1)}] $\Omega_{\mathfrak{g}}(\lambda,\alpha,\beta,\gamma,\bm{\kappa})$ is simple if and only if $\alpha \neq 0\ {\rm or}\ \gamma\neq 0\ {\rm or}\ \bm{\kappa} \neq \bm{0}$ (if $\gamma, \bm{\kappa}$ exist);
\item[{\rm(2)}] $\Omega_{\mathfrak{g}}(\lambda,0,\beta,0,\bm{0})$ has a unique proper submodule $t\Omega_{\mathfrak{g}}(\lambda,0,\beta,0,\bm{0})\cong\Omega_{\mathfrak{g}}(\lambda,1,\beta,0,\bm{0})$, and the quotient $\Omega_{\mathfrak{g}}(\lambda,0,\beta,0,\bm{0})/t\Omega_{\mathfrak{g}}(\lambda,0,\beta,0,\bm{0})$ is a one-dimensional trivial $\mathfrak{g}$-module;
\item[{\rm(3)}] $\Omega_{\mathfrak{g}}(\lambda_1,\alpha_1,\beta_1,\gamma_1,\bm{\kappa}_{1})\cong\Omega_{\mathfrak{g}}(\lambda_2,\alpha_2,\beta_2,\gamma_2,\bm{\kappa}_{2})$ if and only if they have the same parameters, where $\bm{\kappa}_{s}=\{(\kappa_{s,i})_{i\in\Z} \mid \kappa_{s,i}\in\C\},\ s = 1,2$.
\end{itemize}
\end{theo}

\ni{\it Proof.}\ \
The necessity of the statement (1) and the statement (2) can be easily obtained from Lemma~\ref{we-sim&iso}(1) and (2).
Next, we prove the sufficiency of the statement (1).

Assume that $\alpha \neq 0\ {\rm or}\ \gamma\neq 0\ {\rm or}\ \bm{\kappa} \neq \bm{0}$ (if $\gamma, \bm{\kappa}$ exist).
If $\alpha \neq 0$, then by Lemma~\ref{we-sim&iso}(1), $\Omega_{\mathfrak{g}}(\lambda,\alpha,\beta,\gamma,\bm{\kappa})$ is simple when viewed as a $\mathcal{W}(-1)$-module,
and thus it is also a simple $\mathfrak{g}$-module. If $\gamma$ exists and $\gamma\neq 0$, then by Lemma~\ref{hvab-sim&iso}(1),
$\Omega_{\mathfrak{g}}(\lambda,\alpha,\beta,\gamma,\bm{\kappa})$ is simple when viewed as a $\mathfrak{hv}(a,b)$-module, and thus it is also a simple $\mathfrak{g}$-module.

If $\bm{\kappa}$ exists and $\bm{\kappa} \neq \bm{0}$, then there exists at least one $i$ satisfying $\kappa_{i} \neq 0$.
Suppose that $M'$ is a nonzero submodule of $\Omega_{\mathfrak{g}}(\lambda,\alpha,\beta,\gamma,\bm{\kappa})$.
Let $$d_0=\min\{\deg\,f(t)\,|\,f(t)\in M'\}.$$
Since we have assume that $\bm{\kappa}$ exists, we only need to consider the case $b=1$.
If $b = 1, a \neq 0$, by \eqref{action-hvab-1-H} and \eqref{fun-varphi}, we see that
\begin{eqnarray}
\nonumber H_{i,1}\cdot t^{k} &\!\!\!=\!\!\!& \sum_{s=0}^{k}(-1)^s s!\binom{k}{s}\Bigg(\frac{(-1)^{s}}{s!}(a+i)^{s}\kappa_i \\
\nonumber &\!\!\!\!\!\!& +\sum_{\ell=0}^{s-1}\frac{(s-\ell-1)!}{s!}(-1)^{\ell}(a+i)^{\ell}a\lambda^i\beta^{s-\ell}\gamma\Bigg)(t-a-i)^{k-s}\\
\label{Hi1td0-a-nonzero} &\!\!\!=\!\!\!& \kappa_i t^{k} + \mbox{lower degree terms}.
\end{eqnarray}
Similarly, if $b = 1, a=0$, by \eqref{action-hvab-1-H} and \eqref{fun-psi}, we have
\begin{eqnarray}
\nonumber H_{i,1}\cdot t^{k} &\!\!\!=\!\!\!& \sum_{s=0}^{k}(-1)^s s!\binom{k}{s}\frac{(-i)^{s}}{s!}\kappa_i(t-i)^{k-s}\\
\label{Hi1td0-a=0}&\!\!\!=\!\!\!&  \kappa_i t^{k} + \mbox{lower degree terms}.
\end{eqnarray}
Assume that $0\neq f(t)\in M'$ with $\deg\,f(t)=d_0>0$.
Then $\tilde{f}(t):=H_{i,1}\cdot f(t) - \kappa_if(t)\in M'$.
By \eqref{Hi1td0-a-nonzero} and \eqref{Hi1td0-a=0}, we see that $\deg\,\tilde{f}(t)\le d_0-1$.
This contradicts to the definition of $d_0$. Hence, $f(t)$ is a nonzero constant, and thus $1\in M'$, leading to $M'=\Omega_{\mathfrak{g}}(\lambda,\alpha,\beta,\gamma,\bm{\kappa})$.

Finally, we prove the statement (3).
Let $\Omega_s=\Omega_{\mathfrak{g}}(\lambda_s,\alpha_s,\beta_s,\gamma_s,\bm{\kappa}_{s})$, $s=1,2$.
Suppose that $\varphi$ is a module isomorphism from $\Omega_1$ to $\Omega_2$.
Viewing $\Omega_1$ and $\Omega_2$ as $\mathcal{W}(\epsilon)$-modules, by Lemma~\ref{we-sim&iso}(3), we have $\lambda_1=\lambda_2$, $\alpha_1=\alpha_2$, $\beta_1=\beta_2$ and
$\varphi(1)\in\C^*$ (see \cite[Theorems~3.6 and 4.7]{XMDZ2023} for more details).
If $\gamma_1$ and $\gamma_2$ exist, then $b=0$ or $1$. Viewing $\Omega_s$ as $\mathfrak{hv}(a,b)$-modules, by Lemma~\ref{hvab-sim&iso}(3), we have $\gamma_1=\gamma_2$.
If $\bm{\kappa}_{1}$ and $\bm{\kappa}_{2}$ exist, then $b=1$.
For any $i\in\Z$, by \eqref{action-hvab-1-H}, \eqref{fun-varphi} and \eqref{fun-psi}, we always have
$$
\kappa_{1,i}\varphi(1)=\varphi(\kappa_{1,i})=\varphi(H_{i,1}\cdot1)=H_{i,1}\cdot\varphi(1)=\kappa_{2,i}\varphi(1),
$$
which implies $\kappa_{1,i} = \kappa_{2,i}$. Thus, $\bm{\kappa}_{1}=\bm{\kappa}_{2}$.
This completes the proof.
\QED

\vskip15pt

\small \ni{\bf Acknowledgement}
The first author would like to thank Professors S. Tan and Q. Wang for the hospitality during the visit to
Tianyuan Mathematical Center in Southeast China. 
This work was supported by the Fundamental Research Funds for the Central Universities (No.~2024KYJD2006),
and the National Natural Science Foundation of China (No.~12361006).

\end{document}